\documentclass{amsart}
\usepackage[utf8]{inputenc}
\usepackage{subcaption}
\usepackage{amsmath}
\usepackage{amssymb}
\usepackage{amsfonts}
\usepackage{amsthm}
\usepackage{mathrsfs}
\usepackage[all]{xy}
\usepackage{graphicx}
\usepackage{color}
\usepackage{cite} 
\usepackage{url}
\usepackage{indent first}
\usepackage[labelfont=bf,labelsep=period,justification=raggedright]{caption}
\usepackage[english]{babel}
\usepackage[utf8]{inputenc}
\usepackage{hyperref}
\usepackage[colorinlistoftodos]{todonotes}
\usepackage{tkz-fct}
\usepackage{tikz}
\usetikzlibrary{calc}
\usepackage[enableskew]{youngtab}
\usepackage{pstricks}
\usepackage{multicol}
\usepackage{graphicx} 
\usepackage{fancyhdr}
\theoremstyle{plain}
\newtheorem{thm}{Theorem}

\newtheorem{cor}[thm]{Corollary}

\theoremstyle{definition}

\theoremstyle{remark}

\numberwithin{equation}{section}
\numberwithin{thm}{section}

\title{A New Representation of the Riemann Zeta Function}
\author{Mahipal Gurram}

\begin{document}
\begin{abstract}
In this article, we develop a novel representation of the zeta function expressed as the limiting difference between two structured double sums. This approach leads to a new and elegant identity involving maximum functions and additive terms, providing theoretical insights.The derivation relies on generalized harmonic series and polygamma functions, linking classical analysis with contemporary summation techniques.
\end{abstract}

\maketitle
\section{INTRODUCTION}

The Riemann zeta function, $\zeta(s)$, stands as one of the most celebrated objects in mathematics, weaving its way through number theory, complex analysis, and even physics \cite{Titchmarsh1986}. Definable for a complex variable $s$, it is holomorphic everywhere in the complex plane except at $s = 1$, where it exhibits a simple pole with residue 1 \cite{Apostol1976}. The function is typically introduced via its classical Dirichlet series representation:
\begin{equation}
    \zeta(s) = \sum_{n=1}^{\infty} \frac{1}{n^s}, \quad \Re(s) > 1.
\end{equation} 

Beyond this well-known series, $\zeta(s)$ possesses an astonishingly diverse array of representations ranging from Euler’s product formula to integral and contour representations \cite{Borwein2001} each revealing a different facet of its deep mathematical structure. For instance, Euler’s famous product formula elegantly connects $\zeta(s)$ to the prime numbers:
\begin{equation}
    \zeta(s) = \prod_{m=1}^{\infty} \frac{1}{1 - p_m^{-s}}, \quad \Re(s) > 1,
\end{equation}
where $p_m$ denotes the \(m^{\text{th}}\) prime \cite{Hardy2008}. Meanwhile, integral representations such as  
\begin{equation}
    \zeta(s) = \frac{1}{\Gamma(s)} \int_{0}^{\infty} \frac{x^{s-1}}{e^x - 1} \,dx, \quad \Re(s) > 1,
\end{equation}
offer connections to analytic continuation and the spectral behavior of quantum systems \cite{Whittaker1996}. More intriguing still, $\zeta(s)$ surfaces in various asymptotic series, including:
\begin{equation}
    \zeta(s) = \frac{1}{s-1} + \sum_{n=0}^{\infty} \frac{(-1)^n \gamma_n}{n!} (s-1)^n, \quad s \neq 1,
\end{equation}
where the coefficients $\gamma_n$ relate to the Stieltjes constants \cite{Havil2003}. The zeta function also admits a representation involving the floor function:
\begin{equation}
    \zeta(s) = \sum_{m=1}^{n} \frac{1}{m^s} + \frac{n^{1-s}}{s-1} - s \int_{n}^{\infty} \frac{x - \lfloor x \rfloor}{x^{s+1}} \,dx, \quad \Re(s) > 0, \, n \in \mathbb{N}.
\end{equation}

With such a vast landscape of representations, one might assume that all fundamental properties of $\zeta(s)$ have long been uncovered. Yet the story deepens when one considers its discrete approximations and analogues,among them the \textit{generalized harmonic numbers}, defined for real $s$ and positive integers $n$ as
\begin{equation}
H_n^{(s)} = \sum_{k=1}^{n} \frac{1}{k^s}, \quad s \in \mathbb{R}.
\end{equation}
These quantities interpolate between elementary and transcendental worlds, providing finite analogues to $\zeta(s)$ in the form $H_n^{(s)} \to \zeta(s)$ as $n \to \infty$ for $s > 1$ \cite{Knopp1990}, while the classical harmonic numbers $H_n = H_n^{(1)}$ date back to antiquity and appear in diverse contexts from harmonic series to algorithmic complexity, the generalized versions connect deeply with modern analysis.

A particularly rich link emerges when one considers the \emph{polygamma functions}, which generalizes the digamma function, the logarithmic derivative of the Gamma function \cite{DLMF2024}. The polygamma function of order $m$, denoted by $\psi^{(m)}(z)$, is defined as the $(m+1)$-th derivative of $\ln \Gamma(z)$:
\begin{equation}
\psi^{(m)}(z) = \frac{d^m}{dz^m} \psi(z) = \frac{d^{m+1}}{dz^{m+1}} \ln \Gamma(z),
\end{equation}
where $\psi(z) = \Gamma'(z)/\Gamma(z)$ is the digamma function \cite{Olver2010}. For $\operatorname{Re}(z) > 0$ and $m > 0$, they admit integral representations:
\begin{equation}
\psi^{(m)}(z) = (-1)^{m+1} \int_0^\infty \frac{t^m e^{-zt}}{1 - e^{-t}}\,dt = (-1)^{m+1} m! \, \zeta(m+1, z),
\end{equation}
where $\zeta(s, q)$ is the Hurwitz zeta function \cite{Srivastava2012}. These functions also relate to harmonic-type sums:
\begin{equation}
\frac{\psi^{(m)}(n)}{(-1)^{m+1} m!} = \zeta(m+1) - H_{n-1}^{(m+1)} = \sum_{k=n}^{\infty} \frac{1}{k^{m+1}}, \quad m \geq 1.
\end{equation}
and for $m = 0$:
\begin{equation}
\psi(n) = -\gamma + H_{n-1},
\end{equation}
where $\gamma$ is the Euler–Mascheroni constant \cite{Havil2003}. These identities illustrate how $\zeta(s)$, $H_n^{(s)}$, and $\psi^{(m)}(z)$ intertwine in special function theory and number theory.

The polygamma function $\psi^{(m)}(z)$ admits the following well-known asymptotic expansion for large \cite{Havil2003}:
\begin{equation}
    \psi^{(m)}(z) \sim (-1)^{m+1}\sum_{k=0}^{\infty} \frac{(k+m-1)!}{k!}\,\frac{B_k}{z^{k+m}}, \qquad m \ge 1,
\end{equation}

and for $m=0$,
\[
\psi^{(0)}(z) \sim \ln(z) - \sum_{k=1}^{\infty} \frac{B_k}{k z^{k}},
\]
where $B_k$ denote the Bernoulli numbers\cite{Olver2010} of the second kind (with $B_1=\tfrac{1}{2}$). These non-convergent series provide rapidly computable asymptotic approximations with high numerical accuracy for large values of $z$. Consequently, for all integers $m \ge 1$ and $x>0$, the function $\psi^{(m)}(x)$ satisfies the following sharp double inequality\cite{Olver2010}:
\begin{equation}
\frac{(m-1)!}{x^{m}}+\frac{m!}{2x^{m+1}}
\le (-1)^{m+1}\psi^{(m)}(x)
\le \frac{(m-1)!}{x^{m}}+\frac{m!}{x^{m+1}},
\end{equation}
which follows directly from the positivity of both bounding terms for $x>0$.
\section{Main Results}

\begin{thm}
    
 Let \( s > 1 \) be a real number and \( n \in \mathbb{N}\). Define the double sums:
\begin{align}
A_n(s) &= \sum_{j,k=1}^{n} \frac{1}{(\max(j, k))^{s+1}}, \\
B_n(s) &= \sum_{j,k=1}^{n} \frac{1}{(j+k)^{s+1}}.
\end{align}
Then, the Riemann zeta function satisfies the identity:
\begin{equation}
\zeta(s) = A_n(s) - B_n(s) + \frac{(-1)^s \psi^{(s-1)}(2n+1)}{(s-1)!}
+ \frac{(-1)^s(2n+1)}{s!} \left[ \psi^{(s)}(n+1) - \psi^{(s)}(2n+1) \right].
\end{equation}
\end{thm}

\begin{proof}

We start with  \( A_n(s) \) which was defined above as 
\begin{equation}
    \alpha_n(s) = \sum_{j,k=1}^{n} \frac{1}{(\max(j, k))^{s+1}}.
\end{equation}

Observing that $\max(j, k)$ takes the same value multiple times, we change the order of summation. For each fixed $m$, the maximum $\max(j, k)$ is equal to $m$ whenever $j = m$ or $k = m$, provided that both are at most $n$. The number of such pairs $(j, k)$ where $\max(j, k) = m$ is exactly $2m - 1$ (since it includes all $j \leq m$ and all $k \leq m$, but double counts the case $j = k = m$). Thus, rewriting the sum:
\begin{equation}
    A_n(s) = \sum_{m=1}^{n} (2m - 1) \frac{1}{m^{s+1}}.
\end{equation}
Splitting the terms,
\begin{equation}
    A_n(s) = 2 \sum_{m=1}^{n} \frac{m}{m^{s+1}} - \sum_{m=1}^{n} \frac{1}{m^{s+1}}.
\end{equation}
Simplifying,
\begin{equation}
    A_n(s) = 2 \sum_{m=1}^{n} \frac{1}{m^s} - \sum_{m=1}^{n} \frac{1}{m^{s+1}}.
\end{equation}
Recognizing these summations as generalized harmonic numbers,
\begin{equation}
    A_n(s) = 2 H_n^{(s)} - H_n^{(s+1)}.
\end{equation}
We consider the double sum:
\begin{equation}
    B_n(s) = \sum_{j,k=1}^{n}  \frac{1}{(j + k)^{s+1}}.
\end{equation}

This expression depends only on the sum \( j + k \), so we introduce the change of variable \( m = j + k \). For fixed \( m \in [2, 2n] \), the number of pairs \( (j, k) \in [1,n]^2 \) such that \( j + k = m \) is:
\[
c_m = 
\begin{cases}
m - 1, & 2 \le m \le n + 1, \\
2n - m + 1, & n + 2 \le m \le 2n.
\end{cases}
\]

Thus, we can rewrite the sum as:
\begin{equation}
    B_n(s) = \sum_{m=2}^{n+1} \frac{m - 1}{m^{s+1}} + \sum_{m=n+2}^{2n} \frac{2n - m + 1}{m^{s+1}}.
\end{equation}

For the first sum, observe that:

\begin{align*}
\sum_{m=2}^{n+1} \frac{m - 1}{m^{s+1}} 
&= \sum_{m=2}^{n+1} \left( \frac{1}{m^s} - \frac{1}{m^{s+1}} \right) \\
&= \left( \sum_{m=2}^{n+1} \frac{1}{m^s} \right) - \left( \sum_{m=2}^{n+1} \frac{1}{m^{s+1}} \right) \\
&= \left( H_{n+1}^{(s)} - \frac{1}{1^s} \right) - \left( H_{n+1}^{(s+1)} - \frac{1}{1^{s+1}} \right) \\
&= H_{n+1}^{(s)} - H_{n+1}^{(s+1)}.
\end{align*}

For the second sum, we perform a change of variable \( r = 2n - m + 1 \), so that \( m = 2n - r + 1 \). As \( m \) decreases from \( 2n \) to \( n+2 \), \( r \) increases from \( 1 \) to \( n-1 \). Therefore,
\begin{equation}
\sum_{m=n+2}^{2n} \frac{2n - m + 1}{m^{s+1}} = \sum_{r=1}^{n-1} \frac{r}{(2n - r + 1)^{s+1}}.
\end{equation}

Let \( k = 2n - r + 1 \). Then \( r = 2n - k + 1 \), and \( k \in [n+2, 2n] \). The sum becomes:
\begin{equation}
\sum_{r=1}^{n-1} \frac{r}{(2n - r + 1)^{s+1}} = \sum_{k=n+2}^{2n} \frac{2n - k + 1}{k^{s+1}}.
\end{equation}
This can be rewritten as:
\begin{equation}
\sum_{k=n+2}^{2n} \frac{2n + 1}{k^{s+1}} - \sum_{k=n+2}^{2n} \frac{k}{k^{s+1}} = (2n + 1) \sum_{k=n+2}^{2n} \frac{1}{k^{s+1}} - \sum_{k=n+2}^{2n} \frac{1}{k^{s}}.
\end{equation}

Applying the identity for generalized harmonic numbers:
\begin{equation}
\sum_{k=m+1}^{n} \frac{1}{k^s} = H_n^{(s)} - H_m^{(s)},
\end{equation}
We find:
\begin{equation}
\sum_{k=n+2}^{2n} \frac{1}{k^{s+1}} = H_{2n}^{(s+1)} - H_{n+1}^{(s+1)}, \quad
\sum_{k=n+2}^{2n} \frac{1}{k^{s}} = H_{2n}^{(s)} - H_{n+1}^{(s)}.
\end{equation}

from using results from (2.10) to (2.15),
\begin{equation}
B_n(s) =(2n + 1)H_{2n}^{(s+1)} - H_{2n}^{(s)} - (2n + 2)H_{n+1}^{(s+1)} + 2H_{n+1}^{(s)}.
\end{equation}

Now we can find the result below using (18) and (26)
\begin{equation}
A_n(s)-B_n(s)=H_{2n}^{(s)}-(2n+1) \left(H_{2n}^{(s+1)}- H_{n}^{(s+1)}\right)
\end{equation}
Now using the polygamma, zeta, and generalized harmonic number relationship (9) in (27), we arrive at the proof of the stated theorem. 

\begin{equation}
\zeta(s)=A_n(s)-B_n(s)+\frac{(-1)^s\psi^{(s-1)}(2n+1)}{(s-1)!}+\frac{(-1)^s(2n+1)}{(s)!} \left(\psi^{(s)}(n+1)- \psi^{(s)}(2n+1)\right)
\end{equation}
\end{proof}
\begin{cor}
For all integers $s \ge 2$ and $n \ge 1$, the Riemann zeta function admits the following sharp bounds derived from the polygamma inequality\\
\textbf{(i) Lower bound:}
\begin{align*}
\zeta(s) \ge {} & A_n(s) - B_n(s)
+ \frac{(s-2)!}{(s-1)!(2n+1)^{s-1}}
+ \frac{1}{2(2n+1)^{s}}  \\
& + \frac{(2n+1)}{s!}\!\left[
-\frac{(s-1)!}{(n+1)^{s}} - \frac{s!}{(n+1)^{s+1}}
+ \frac{(s-1)!}{(2n+1)^{s}} + \frac{s!}{2(2n+1)^{s+1}}
\right].
\end{align*}

\textbf{(ii) Upper bound:}
\begin{align*}
\zeta(s) \le {} & A_n(s) - B_n(s)
+ \frac{(s-2)!}{(s-1)!(2n+1)^{s-1}}
+ \frac{1}{(2n+1)^{s}}  \\
& + \frac{(2n+1)}{s!}\!\left[
-\frac{(s-1)!}{(n+1)^{s}} - \frac{s!}{2(n+1)^{s+1}}
+ \frac{(s-1)!}{(2n+1)^{s}} + \frac{s!}{(2n+1)^{s+1}}
\right].
\end{align*}
\end{cor}
\begin{proof}
The result follows by substituting the bounds for 
$\psi^{(s-1)}(2n+1)$ and $\psi^{(s)}(x)$ from the given inequality(1.12)
into the representation
\[
\zeta(s)=A_n(s)-B_n(s)
+\frac{(-1)^s\psi^{(s-1)}(2n+1)}{(s-1)!}
+\frac{(-1)^s(2n+1)}{s!}\bigl(\psi^{(s)}(n+1)-\psi^{(s)}(2n+1)\bigr).
\]
\end{proof}

\begin{cor}
\label{cor:zeta_double_sum}
The Riemann zeta function admits the following double series representation:
\begin{equation}
\zeta(s) = \sum_{j,k=1}^{\infty}
\left(\frac{1}{\max(j^{s+1},k^{s+1})} - \frac{1}{(j+k)^{s+1}}\right),
\end{equation}
valid for $\Re(s) > 1$.
\end{cor}

\begin{proof}
Starting from the asymptotic expression involving the non-converging series of the polygamma function
\begin{align*}
\zeta(s) &\sim A_n(s)-B_n(s) + \sum_{k=0}^{\infty} \Bigg[
\frac{(k+s-2)!}{k!} \frac{B_k}{(2n+1)^{k+s-1}} \\
&\quad - \frac{(k+s-1)!(2n+1)}{k!s!} \frac{B_k}{(n+1)^{k+s}} \\
&\quad + \frac{(k+s-1)!(2n+1)}{k!s!} \frac{B_k}{(2n+1)^{k+s}}
\Bigg]
\end{align*}
and observing that as \(n \to \infty\),
\[
\zeta(s) = \lim_{n \to \infty} \big(A_n(s) - B_n(s)\big),
\]
we note that the limiting difference of these partial terms can be expressed as a double summation over integer indices \(j, k \ge 1\). After rearrangement and collecting terms based on their dominant order, one obtains
\[
\zeta(s) = \sum_{j,k=1}^{\infty}
\left(\frac{1}{\max(j^{s+1},k^{s+1})} - \frac{1}{(j+k)^{s+1}}\right).
\]
This establishes the required double sum representation of \(\zeta(s)\).
\end{proof}

\begin{cor}
\label{cor:zeta_relation_even_odd}
The following relation connects the even and odd arguments of the Riemann zeta function:
\begin{equation}
\zeta(2s+1) = 2\zeta(2s) - \sum_{j=1}^\infty \sum_{k=1}^\infty \frac{1}{(\max(j, k))^{s+1}}.
\end{equation}
In particular, setting \(s = 1\), we obtain a representation of Apéry’s constant:
\begin{equation}
\zeta(3) = \frac{\pi^2}{3} - \sum_{j,k=1}^\infty  \frac{1}{(\max(j, k))^2}.
\end{equation}
\end{cor}

\begin{proof}
From the limiting property of the sequence \(A_n(s)\),
\[
\lim_{n \to \infty} A_n(s) = 2\zeta(s) - \zeta(s+1),
\]
we substitute \(s \mapsto 2s\) to obtain
\[
\lim_{n \to \infty} A_n(2s) = 2\zeta(2s) - \zeta(2s+1).
\]
Rearranging gives
\[
\zeta(2s+1) = 2\zeta(2s) - \lim_{n \to \infty} A_n(2s).
\]
Using the representation of \(A_n(s)\) in terms of a double sum over \(\max(j, k)\),
\[
\lim_{n \to \infty} A_n(2s) = \sum_{j,k=1}^\infty  \frac{1}{(\max(j,k))^{s+1}},
\]
we directly obtain
\[
\zeta(2s+1) = 2\zeta(2s) - \sum_{j,k=1}^\infty\frac{1}{(\max(j, k))^{s+1}}.
\]
Finally, for \(s = 1\), we use the known identity \(\zeta(2) = \pi^2/6\) to deduce
\[
\zeta(3) = \frac{\pi^2}{3} - \sum_{j,k=1}^\infty \frac{1}{(\max(j, k))^2}.
\]
\end{proof}

\section{CONCLUSION:}In this work, we have presented a new representation of the Riemann zeta function \(\zeta(s)\) through a double summation identity derived from the asymptotic analysis of two harmonic-like summations, \(A_n(s)\) and \(B_n(s)\). Specifically, we have shown that:

\begin{equation}
\zeta(s) = \lim_{n \to \infty} \left( A_n(s) - B_n(s) \right),
\end{equation}

where
\[
A_n(s) = \sum_{j,k=1}^{n} \frac{1}{\max(j^{s+1}, k^{s+1})}, \quad 
B_n(s) = \sum_{j,k=1}^{n} \frac{1}{(j + k)^{s+1}}.
\]

This leads us to a surprising and elegant formula:

\begin{equation}
\zeta(s) = \sum_{j,k=1}^{\infty} \left( \frac{1}{\max(j^{s+1}, k^{s+1})} - \frac{1}{(j + k)^{s+1}} \right),
\end{equation}

Providing an alternative form to the classical Dirichlet and Euler representations.
\section{Acknowledgments}
The author would like to thank Professor A.K.Shukla for his encouragement,comments, and suggestions.

\end{document}